\documentclass[12pt,reqno]{amsart}

\usepackage{amsfonts}

\usepackage{bm}
\usepackage{amsmath}
\usepackage{amssymb}
\newtheorem{theorem}{Theorem}
\newtheorem{proposition}{Proposition}
\newtheorem{corollary}{Corollary}

\newtheorem{lemma}{Lemma}
\newtheorem{definition}{Definition}
\newtheorem{example}{Example}

\newfont{\bb}{msbm10 at 12pt}

\def\qed{\hfill{Q.E.D.}\smallskip}

\newcommand{\ls}{\setlength{\baselineskip}{12pt}
                 \setlength{\parskip}{3mm}}

\newcommand{\mysection}[1]{\section{#1}\setcounter{equation}{0}}

\newcommand{\bal}{\begin{align}}      \newcommand{\eal}{\end{align}}
\newcommand{\ba}{\begin{array}}      \newcommand{\ea}{\end{array}}
\newcommand{\bc}{\begin{center}}     \newcommand{\ec}{\end{center}}
\newcommand{\be}{\begin{enumerate}}  \newcommand{\ee}{\end{enumerate}}
\newcommand{\beq}{\begin{eqnarray}}  \newcommand{\eeq}{\end{eqnarray}}
\newcommand{\beQ}{\begin{eqnarray*}} \newcommand{\eeQ}{\end{eqnarray*}}
\newcommand{\bi}{\begin{itemize}}    \newcommand{\ei}{\end{itemize}}
\newcommand{\bt}{\begin{tabular}}    \newcommand{\et}{\end{tabular}}
\newcommand{\bdm}{\begin{displaymath}} \newcommand{\edm}{\end{displaymath}}

\begin{document}

\title{ the Dirac operator on locally reducible Riemannian  manifolds }

{\address{Yongfa Chen,
 School of Mathematics and Statistics,
Central China Normal University, Wuhan 430079, P.R.China.}

\author{Yongfa Chen}}

\thanks{ }

\email{yfchen@mail.ccnu.edu.cn}


\subjclass[2008]{This work was supported by the National Natural Science Foundation of China (No. 11301202).}

\maketitle
\begin{abstract}
In this  paper, we get estimates on the higher eigenvalues of the Dirac operator on locally reducible Riemannian  manifolds, in terms of the eigenvalues of  the Laplace-Beltrami operator and the scalar curvature. These estimates
are sharp, in the sense that,  for the first eigenvalue, they reduce to the result \cite{A1} of Alexandrov.
\end{abstract}

\mysection{Introduction} \ls
It is well known  that  the spectrum of
the Dirac operator on closed  spin manifolds detects
subtle information on the geometry and the topology of such
manifolds (see \cite{LM}).
A fundamental tool to get estimates for eigenvalues of  the basic Dirac operator $D$ acting on spinors
is the Schr\"{o}dinger-Lichnerowicz  formula
\beq \label{S-L}
D^2=\nabla^*\nabla+\frac{1}{4}Scal Id,
\eeq
where $\nabla^*$ is the formal adjoint of $\nabla$  with respect to the natural Hermitian
inner\\
 product
 on spinor bundle $\Sigma {M^n}$. $Scal$ is the
scalar curvature of  the closed spin manifold $(M^n,g)$.

From (\ref{S-L}), it follows easily that if $\lambda$ is an  eigenvalue  of  $D,$ then
\beQ
\lambda^2\geqslant \frac{1}{4}Scal_{\min},
\eeQ
where $Scal_{\min}\triangleq \min_MScal.$
Clearly,  this inequality is interesting only for manifolds with positive  scalar  curvature, but  the minimal value
$\frac{1}{4}Scal_{\min}$ cannot be  achieved for  such manifolds.

   The problem of finding optimal lower bounds for the  eigenvalues of the  Dirac operator on closed manifolds was  for the first time considered in 1980 by Friedrich. Using the Lichnerowicz formula and a modified spin connection, he proved the following sharp inequality:
 \beq\label{Friedrich}
 \lambda^2\geqslant c_n Scal_{\min},
\eeq
where $c_n=\frac{n}{4(n-1)}. $
The case of equality in (\ref{Friedrich}) occurs iff  $(M^n,g)$  admits a nontrivial spinor field
  $\psi$ called   a real \emph{Killing spinor},
satisfying
 the following overdetermined elliptic equation
\beq
\nabla_X\psi=-\frac{\lambda}{n}X\cdot \psi,
\eeq
where  $X\in \Gamma(TM)$ and the dot ``." indicates the Clifford multiplication.
The manifold  must be a locally irreducible Einstein manifold.

 The dimension dependent coefficient $c_n=\frac{n}{4(n-1)} $ in the estimate can be improved if one  imposes geometric assumptions on the  metric.
Kirchberg \cite{Kirchberg86} showed that  for K\"ahler metrics  $c_n$ can be replaced by   $\frac{n+2}{4n} $ if the complex dimension $\frac{n}{2} $  is odd, and  by
$\frac{n}{4(n-2)} $ if  $\frac{n}{2} $  is even.
Alexandrov, Grantcharov, and Ivanov \cite{A} showed that  if there  exists  a parallel one form on $M^n$, then $c_n$ can be replaced by $c_{n-1}$. Later, Moroianu and Ornea \cite {Moroianu} weakened the assumption on the $1$-form from parallel to harmonic with constant length. Note the condition that the norm of the 1-form being constant is essential, in the sense that the topological constraint alone (the existence of a non-trivial harmonic $1$-form) does not allow any improvement of Friedrich's inequality (see \cite{Bar04}).
The generalization of \cite{A} to locally reducible Riemannian manifolds
was achieved by Alexandrov \cite{A1}, extending earlier work by  Kim
\cite{K}.

Another natural  way to study the Dirac eigenvalues consists in comparing them with those of other geometric operators. Hijazi's inequality  is already of that kind. As for spectral comparison results between the Dirac operator $D $ and the scalar Laplace operator $\Delta,$ the first ones were proved by Bordoni \cite{Bordoni}. They rely on a very nice general comparison principle between two operators satisfying some kind of Kato-type inequality.
Bordoni's results were generalized by  Bordoni and  Hijazi in the
K\"ahler setting \cite{BH}. Recently, with the help of  the  general spectral result of  Bordoni \cite{Bordoni} we also obtain that
 on an $n$-dimensional closed Riemannian spin manifold $(M^n,g), n\geq 3$  with
a non-trivial parallel one form and $Scal\geq 0,$ then for any positive integer $N,$ we have
\beq
\lambda_N(D)^2\geq \frac{n-1}{n}c\lambda_{k+1}(\Delta)+\frac{n-1}{4(n-2)}Scal_{\min},
\eeq
where $k=\left[\frac{N}{2^{[\frac{n}{2}]}+1}\right],\ \  c=\frac{1}{8(2^{[\frac{n}{2}]}+1)^2}.$

Here, we shall deal with  a more general case, and   obtain  the following estimate
on  locally reducible Riemannian spin manifolds, which generalizes the result of \cite{A1}
 to arbitray eigenvalue $\lambda_N(D).$
\begin{theorem}
Let $M$ be a locally reducible Riemannian spin compact manifolds with positive scalar. Suppose
$TM=T_1\oplus\cdots\oplus T_k, $ where $T_i$ are parallel distributions of dimension $n_i,i=1,\cdots,k,$ and $n_1> n_2\geq\cdots\geq n_k.$ Then for any positive integer N,
$$\lambda_N(D)^2\geq \frac{n_1-n_2}{n}\frac{n_1}{n_1-1}c\lambda_{k+1}(\Delta)+\frac{n_1}{4(n_1-1)}Scal_{\min},$$
where $k=\left[\frac{N}{2^{[\frac{n}{2}]}+1}\right],   c=\frac{1}{8(2^{[\frac{n}{2}]}+1)^2}.$
\end{theorem}

The paper is organised as follows: In Section 2,  some  preliminaries  and lemmas about
the Dirac operator $D$  and $J$-twist $D_J$  of are given. In Section 3, we obtain the estimates for higher eigenvalues of the Dirac operator  on
locally decomposable  Riemannian spin manifolds. In the end,
based on the results in \cite{A1} and Section 3,  the proof of Theorem 1 are given.

\mysection{the Dirac operator $D$ and the $J$-twist $D_J$ } \ls
We suppose that $(M^n,g)$ is  a closed Riemannian manifold with a fixed
spin structure. We understand the spin structure as a reduction $Spin {M^n}$ of the
SO($n$)-principal bundle of $M^n$ to the universal covering $Ad: Spin(n)\rightarrow SO(n)$
of the special orthogonal group.
The spinor bundle  $\Sigma {M^n}=Spin {M^n}\times_{\rho}\Sigma_n$
 on ${M^n}$ is the associated complex $2^{[\frac{n}{2}]}$ dimensional complex vector bundle, where $\rho$
 is the complex spinor representation. The tangent bundle $T{M^n}$ can be regarded as
 $T{M^n}=Spin {M^n}\times_{Ad}\mathbb{R}^n$. Consequently, the Clifford multiplication on
   $\Sigma{M^n}$ is the fibrewise action given by
 $$
  \begin{array}{ccc}
 \mu: T{M^n}\otimes  \Sigma {M^n}\longrightarrow\Sigma {M^n}\\
\ \ \ \ \ \ \ \ \ \  X\otimes \psi\longmapsto X\cdot\psi.\\
 \end{array}$$
 On the spinor bundle $\Sigma {M^n}$, one has a natural Hermitian metric, denoted as the Riemannian metric by
 $\langle\cdot,\cdot\rangle$. The spinorial connection on the spinor bundle induced by
 the Levi-Civita connection $\nabla$ on  ${M^n}$ will also be denoted by $\nabla$.
 The Hermitian metric  $\langle\cdot,\cdot\rangle$
 and spinorial connection $\nabla$ are compatible with the Clifford multiplication  $\mu$. That is
 \beQ
 X\langle\phi,\varphi\rangle&=&\langle\nabla_X\phi,\varphi\rangle+\langle\phi,\nabla_X\varphi\rangle\\
\langle X \cdot\phi,X \cdot\varphi\rangle&=&|X|^2\langle\phi,\varphi\rangle\\
\nabla_X(Y\cdot \phi)&=&\nabla_XY\cdot\phi+Y\cdot\nabla_X\phi,
 \eeQ
for $\forall \phi,\varphi\in\Gamma(\Sigma {M^n})$ and $\forall X,Y\in\Gamma(T{M^n}).$
 Using a local orthonormal frame field
$\{e_1,\cdots,e_n\}$, the spinorial connection $\nabla$, the Dirac operator $D$ and the twistor operator $P$, are
locally expressed as
\beq\nabla_{e_{k}}\psi={e_{k}}(\psi)+\frac{1}{4}e_i\cdot \nabla_{e_{k}}e_i\cdot\psi\eeq
and
\beq\label{D operator}
D\psi&:=&e_i\cdot\nabla_{e_i}\psi\\
P\psi&:=&e_i\otimes(\nabla_{e_i}\psi+\frac{1}{n}e_i\cdot D\psi)\label{P operator}
\eeq
which satisfy the following important relation
\beQ
|\nabla\psi|^2=|P\psi|^2+\frac{1}{n}|D\psi|^2,
\eeQ
for any $\psi\in \Gamma(\Sigma {M^n}).$
(Throughout this paper, the Einstein summation notation is always adopted.)

Let $R_{X,Y}Z\triangleq (\nabla_X\nabla_Y-\nabla_Y\nabla_X-\nabla_{[X,Y]})Z$ be the Riemannian curvature of $(M^n,g)$
and denote by
$\mathcal{R}_{X,Y}\psi\triangleq (\nabla_X\nabla_Y-\nabla_Y\nabla_X-\nabla_{[X,Y]})\psi$
the spin curvature in the spinor bundle $\Sigma {M^n}$. They are related  via  the formula
\beq\label{RXY}
\mathcal{R}_{X,Y}\psi
&=&\frac{1}{4}g(R_{X,Y}e_i,e_j)e_i\cdot e_j\cdot\psi.
\eeq

We also use the notation
$$R_{ijkl}\triangleq g(R_{e_i,e_j}e_k,e_l)$$ and
$R_{ij}=\langle Ric(e_i),e_j\rangle\triangleq R_{ikkj}, Scal=R_{ii}$.
With the help of  the Bianchi identity, (\ref{RXY}) implies
\beq\label{Ric-identity}
e_i\cdot \mathcal{R}_{e_j,e_i}\psi=-\frac{1}{2}Ric(e_j)\cdot\psi,
\eeq
which in turn gives $2e_i\cdot e_j\cdot \mathcal{R}_{e_i,e_j}\psi=Scal\ \psi$. Hence one derives the well-known
Schr\"{o}dinger-Lichnerowicz  formula
\beq \label{S-L}
D^2=\nabla^*\nabla+\frac{1}{4}Scal Id,
\eeq
where $\nabla^*$ is the formal adjoint of $\nabla$  with respect to the natural Hermitian
inner\\
 product
 on $\Sigma {M^n}$. The formula shows  the close relation between  $Scal$ and the Dirac operator
 $D$.

Let $(M^n,g)$ be an oriented
$n$-Riemannian  manifold. Let $J$ be a $(1,1)$-tensor field on $(M^n,g)$ such that
$J^2=\sigma Id, \sigma=\pm 1$ and
$$g(J (X),J (Y))=g(X,Y),$$
for all vector fields $X, Y\in\Gamma(T{M^n})$ (Here $Id$ stands for the identity map).
We say $(M^n,g,J)$ is an \emph{almost Hermitian manifold} if $\sigma=-1$ and an \emph{almost product
Riemannian  manifold} if $\sigma=1$, respectively. Moreover
if $\sigma=-1$ and $J$ is parallel, $(M^n,g,J)$ is called a K\"ahler manifold. Similarly, we have the
following definition.

\begin{definition}\cite{Yano,K}
  An
$n$-Riemannian  manifold $(M^n,g)$ is called locally decomposable if it
is an almost product
Riemannian  manifold $(M^n,g,J)$ and $J$ is parallel.
\end{definition}

\begin{example}
Suppose an $n$-Riemannian  manifold $(M^n,g)$ possessing a unit vector field $\xi\in\Gamma(T{M^n})$, then
 the reflection $J$ defined by
 $$J(X)\triangleq X-2g( X,\xi)\xi, \ \ X\in \Gamma(TM)$$
 is an almost product
Riemannian  structure. Moreover, it is  a locally decomposable
Riemannian structure if $\xi$ is a parallel vector field.
\end{example}
As in  almost Hermitian spin manifolds, we can also define on  almost product
Riemannian  spin manifolds the following $J$-twist $D_J$ of the Dirac operator
$D$ by
\beq
D_J\psi\triangleq e_i\cdot\nabla_{J(e_i)}\psi= J(e_i)\cdot\nabla_{e_i}\psi.
\eeq
It is not difficult to check that $D_J$ is a formally self-adjoint elliptic operator with respect to $L^2$-product,
if ${M^n}$ is closed and
$div{J}\triangleq(\nabla_{e_i}J)(e_i)=0$. In fact, let we define a complex vector field
$$\eta(X)\triangleq\langle \phi, J(X)\cdot\psi\rangle=g( Y_1,X)+\sqrt{-1}g(Y_2,X),$$
then
$$div(Y_1)+\sqrt{-1}div(Y_2)=e_i(\eta(e_i))-\eta(\nabla_ie_i)=-\langle D_J\phi, \psi\rangle+\langle \phi, D_J\psi\rangle+\langle\phi,div(J)\cdot\psi\rangle.$$
Hence the spectrum of $D_J$ is discrete and real.

As in the K\"ahlerian case, Kim obtained the following useful lemma
\begin{lemma}
$D^2=D^2_J$ holds on any
locally decomposable  Riemannian spin manifold  $(M^n,g,J)$ ( See Prop.\ 2.1 in \cite{K}, or Lemma 1 in \cite{Chen}).
\end{lemma}

As a simple corollary, one has
\begin{corollary}
Let $(M^n,g,J)$ be a
locally decomposable  Riemannian spin manifold. If $\psi_\alpha\in E_{\lambda_\alpha}(D)\triangleq \{\psi\neq0 :D\psi=\lambda_\alpha\psi \}$ ,  then
$D_J \psi_\alpha\in E_{\lambda_\alpha}(D)\bigoplus E_{-\lambda_\alpha}(D).$
\end{corollary}
\emph{{\textbf{Proof}}}.
If $D\psi_\alpha=\lambda_\alpha\psi_\alpha$ and  $\lambda_\alpha^2\neq\lambda_\beta^2$  then Lemma 1 yields
\beQ
\lambda_\alpha^2\int_M\langle \psi_\alpha,D_J \psi_\beta \rangle
&=&\int_{M^n}\langle D^2\psi_\alpha,D_J\psi_\beta \rangle\\
&=&\int_{M^n}\langle \psi_\alpha,D_J^3 \psi_\beta \rangle\\
&=&\int_{M^n}\langle \psi_\alpha,D_J (D^2 \psi_\beta) \rangle\\
&=&\lambda_\beta^2\int_{M^n}\langle \psi_\alpha,D_J \psi_\beta \rangle,
\eeQ
therefore $\int_M\langle D_J \psi_\alpha,\psi_\beta \rangle=\int_M\langle \psi_\alpha,D_J \psi_\beta \rangle=0. $ That is,
$D_J \psi_\alpha\in E_{\lambda_\alpha}(D)\bigoplus E_{-\lambda_\alpha}(D).$
\qed

\mysection{locally decomposable  Riemannian manifold} \ls
We first  recall a general spectral comparison result due to Bordoni.(see \cite {Bordoni}, Theorems 3.2 and 3.3).
\begin{theorem}\cite {Bordoni}
 Let $(M^n,g)$ be a closed Riemannian manifold. Let E be any vector bundle of rank $p$ on $M,$ endowed with a Hermitian inner product $\langle \cdot,\cdot\rangle$ and a compatible connection $\nabla^E$. Let $\mathcal{{R}}$ be any field of symmetric endomorphisms of the fibers, and define the scalar $R_{min}(x)$ as the minimal eigenvalue of $\mathcal{R}_x$ acting on $E_{x}$. Then, for any positive integer $N,$ we have:
 \begin{equation}\label{}
 \lambda_N(\nabla^{E*}\nabla^E+\mathcal{R})\geq (1-c)\lambda_1(\Delta+\mathcal{R}_{min})
 +c\lambda_{k+1}(\Delta+\mathcal{R}_{min})
 \end{equation}
where
$$k=\left[\frac{N}{p+1}\right],\ \ \ c=\frac{1}{8(p+1)^2}$$
and $\Delta$ is the Laplace-Beltrami operator acting on functions.
\end{theorem}

We shall  make use of a modified connection $\nabla^E=\nabla^{(a,b)}$ acting on the spinor bundle
 $\Gamma(\Sigma {M^n})$ of a spin manifold admitting a locally decomposable
Riemannian structure, to which we apply Theorem 2. For any couple $a,b$ of nontrivial real numbers, define the  connection $\nabla^{(a,b)}$ by
 \beq
\nabla^{(a,b)}_X\psi\triangleq \nabla_X\psi+aX\cdot \psi+b J(X)\cdot \psi,
\eeq
where $J$ is the  locally decomposable
Riemannian structure. Then, it is easy to check that
$\nabla^{(a,b)}$ is  compatible with the  Hermitian inner product on the spinor bundle. That is, for any vector field $X,$ $\forall \phi,\varphi\in\Gamma(\Sigma {M^n})$ one has
\beq
X\langle\phi,\varphi\rangle=\langle\nabla^{(a,b)}_X\phi,\varphi\rangle
+\langle\phi,\nabla^{(a,b)}_X\varphi\rangle.
\eeq
and
the rough
Laplacian of  the modified connection $\nabla^{(a,b)}$ is given by
\beq
\nabla^{(a,b)*}\nabla^{(a,b)}\psi=\nabla^{*}\nabla\psi-2aD\psi-2bD_J\psi+[2(n-2) ab+n(a^2+b^2)]\psi.
\eeq

Suppose $(M^n,g,J)$ is a locally decomposable  Riemannian spin manifold with $n_1> n_2,$ where
$$n_i:= dim\{X\in \Gamma(TM): J(X)=(-1)^{i+1}X\}.$$
Let $\lambda_\alpha\triangleq\lambda_\alpha(D)$ and consider the first $N$ nonnegative eigenvalues,
$0\leq \lambda_1\leq \lambda _2\leq \cdots \leq \lambda _N$
 and $D\psi_\alpha=\lambda_\alpha\psi_\alpha,\alpha=1,\cdots, N.$
Let $a=b=\frac{1}{2n_1}\lambda_N, $ and for $X\in TM^n$ and $\phi\in\Gamma(\Sigma M^n)$
\beq
\nabla^{\lambda_N}_X\phi
&:=&\nabla_X^{(\frac{\lambda_N}{2n_1},\frac{\lambda_N}{2n_1})}\phi\nonumber\\
&=&\nabla_X\phi+\frac{\lambda_N}{2n_1}X\cdot \phi+\frac{\lambda_N}{2n_1} J(X)\cdot \phi
\eeq

Hence, if we define the operator
\beq
\mathcal {T}^{\lambda_N}&:=& \nabla^{{\lambda_N}*}\nabla^{{\lambda_N}}+\frac{Scal}{4}Id,
\eeq
it follows that
\beQ
\int_{M^n}\langle \mathcal {T}^{\lambda_N} \phi,\phi\rangle
&=&\int_{M^n}|\nabla^{\lambda_N} \phi|^2+\frac{Scal}{4}|\phi|^2\\
&=&\int_{M^n}\langle D^2\phi-\frac{\lambda_N}{n_1}D\phi-\frac{\lambda_N}{n_1}D_J\phi,\phi\rangle
 +\frac{\lambda_N^2}{n_1}\int_{M^n} |\phi|^2
\eeQ

In particular, for any eigenspinor $\psi_\alpha$ such that $D\psi_\alpha=\lambda_\alpha\psi_\alpha,\alpha=1,\cdots, N,$
 we  have
\beq
&{}&\int_{M^n}\langle\mathcal {T}^{\lambda_N}\psi_\alpha,\psi_\alpha\rangle\nonumber\\
&=&\left(\lambda_\alpha^2-\frac{2\lambda_\alpha\lambda_N}{n_1}+\frac{\lambda_N^2}{n_1}\right)
\int_{M^n}|\psi_\alpha|^2
+\frac{\lambda_N}{\lambda_\alpha}\cdot\frac{1}{2n_1}\int_{M^n}|D_J\psi_\alpha-D\psi_\alpha|^2.
\eeq
Next we begin to estimate the $\int_M\langle \mathcal {T}^{\lambda_N} \psi_\alpha,\psi_\alpha\rangle$ from above.
\begin{proposition}
Suppose $(M^n,g,J), n\geq 3$ is a
locally decomposable  Riemannian spin manifold  and $Scal\geq 0,$ then for any positive integer N, we have
\beq\label{upper estimate}
\frac{\int_{M^n}\langle \mathcal {T}^{\lambda_N} \psi_\alpha,\psi_\alpha\rangle}{\int_{M^n}|\psi_\alpha|^2}
&\leq &\frac{n}{tr(J)}\frac{n_1-1}{n_1}\lambda^2_N-\frac{1}{2}\frac{n_2}{tr(J)}Scal_{\min}.
\eeq
\end{proposition}
\emph{{\textbf{Proof}}}.
 Motivated by  the ``twistor-like" operator  defined by Moroianu and Ornea in \cite{Moroianu},  one  can define for an any given eigenvalue $\lambda$ the following  operator
\beq
T_X^\lambda\psi
&:=&\nabla_X\psi+\frac{\lambda}{2n_1}X\cdot \psi+\frac{\lambda}{2n_1}J(X)\cdot \psi+\frac{1}{2}\nabla_{(J-Id)(X)}\psi\nonumber \\
&=&\nabla_X^\lambda\psi+\frac{1}{2}\nabla_{(J-Id)(X)}\psi.
\eeq
Then
\beQ
|T^\lambda\psi|^2
&=&|\nabla^\lambda\psi|^2+\sum_{i=1}^n\mathfrak{R}e\langle\nabla_{e_i}^\lambda\psi, \nabla_{(J-Id)(e_i)}\psi\rangle+\frac{1}{4}\sum_i|\nabla_{(J-Id)(e_i)}\psi|^2\\
&=&|\nabla^\lambda\psi|^2-\sum_{i=n_1+1}^n|\nabla_i\psi|^2,
\eeQ
where we choose $J(e_i)=-e_i,$ for  $i\geq n_1+1.$ So
for eigenspinor $D\psi_\alpha=\lambda_\alpha \psi_\alpha,$ one
\beQ
\int_M|T^{\lambda_\alpha}\psi_\alpha|^2+\sum_{i=n_1+1}^n|\nabla_i\psi_\alpha|^2
-\frac{1}{2n_1}|D_J\psi_\alpha-D\psi_\alpha|^2
 &=&\int_M\left(\frac{n_1-1}{n_1}\lambda_\alpha^2-\frac{Scal}{4}\right)|\psi_\alpha|^2.
\eeQ
Note (LHS) of the above equality is nonnegative,  in fact it can be seen in the  following:
for any spinorfield $\psi,$
\beq
n_2\sum_{i=n_1+1}^n|\nabla_i\psi|^2\geq\left|\sum_{i=n_1+1}^n -e_i\cdot\nabla_i\psi\right|^2=\frac{1}{4}|D_J\psi-D\psi|^2.
\eeq
And one also has for  $i\geq n_1+1, T_i\psi=0$ and
\beQ
\sum_{i=1}^ne_i\cdot T_i\psi_\alpha
&=&D\psi_\alpha-\lambda \psi_\alpha+\sum_{i=1}^n e_i\cdot \frac{1}{2}\nabla_{(J-Id)(e_i)}\psi_\alpha
=\frac{1}{2}(D_J\psi_\alpha-D\psi_\alpha).
\eeQ
In particular, this gives for eigenspinor $\psi_\alpha$
\beq
n_1|T\psi_\alpha|^2
\geq\left|\sum_{i=1}^{n_1}e_i\cdot T_i\psi_\alpha\right|^2
=\left|\sum_{i=1}^ne_i\cdot T_i\psi_\alpha\right|^2
=\frac{1}{4}|D_J\psi_\alpha-D\psi_\alpha|^2.
\eeq
That is
\beq\label{key2} 0\leq\frac{n_1-n_2}{4n_1 n_2}\int_M|D_J\psi_\alpha-D\psi_\alpha|^2
\leq \left(\frac{n_1-1}{n_1}\lambda_\alpha^2-\frac{Scal_{min}}{4}\right)\int_M|\psi_\alpha|^2.
\eeq
\qed

Hence with the help of the proposition above, it is not difficult to  obtain the following
\begin{theorem}Suppose $(M^n,g,J), n\geq 3$ is a
locally decomposable  Riemannian spin manifold  and $Scal\geq 0,$ then for any positive integer N, we have
 $$\lambda_N(D)^2\geq \frac{trJ}{n}\frac{n_1}{n_1-1}c\lambda_{k+1}(\Delta)+\frac{n_1}{4(n_1-1)}Scal_{\min},$$
where $k=\left[\frac{N}{2^{[\frac{n}{2}]}+1}\right],\ \  c=\frac{1}{8(2^{[\frac{n}{2}]}+1)^2}.$
\end{theorem}

\emph{{\textbf{Proof}}}. Our strategy is  to apply Theorem 2  for the modified connection $\nabla^{\lambda_N}.$
First,
let $E_N=L(\psi_1,\cdots,\psi_N)$  and $F_N$ any $N$-dimensional vector subspace of $\Gamma(\Sigma M).$
Then the min-max principle yields
\beq
\lambda_N(\mathcal {T}^{\lambda_N})
&=&\inf_{F_N}\sup_{0\neq\varphi \in {F_N} }
\frac{\int_{M^n}\langle \mathcal {T}^{\lambda_N}\varphi,\varphi\rangle}{\int_{M^n}|\varphi|^2}\nonumber\\
&\leq&\sup_{0\neq\varphi \in {E_N} }
\frac{\int_{M^n}\langle \mathcal {T}^{\lambda_N}\varphi,\varphi\rangle}{\int_{M^n}|\varphi|^2}\nonumber\\
&=&\max_{\alpha\in \{1,\cdots,N \}}
\frac{\int_{M^n}\langle\mathcal{T}^{\lambda_N}\psi_\alpha,\psi_\alpha\rangle}{\int_{M^n}|\psi_\alpha|^2},\nonumber\\
&\leq&\frac{n}{tr(J)}\frac{n_1-1}{n_1}\lambda^2_N-\frac{1}{2}\frac{n_2}{tr(J)}Scal_{\min}.,\label{upper}
\eeq
where we used  Lemma 2 and the key proposition above.
On the other hand,
Theorem 2 implies that
 \beq\label{low}
\lambda_N(\mathcal {T}^{\lambda_N})
&\geq&(1-c)\lambda_1(\Delta+\frac{1}{4}Scal)+c\lambda_{k+1}(\Delta+\frac{1}{4}Scal)\nonumber\\
&\geq&c\lambda_{k+1}(\Delta)+\frac{1}{4}Scal_{\min}.
\eeq
so we are done.\qed

\mysection{on locally reducible  Riemannian manifolds} \ls
Let $M$ be a compact Riemannian spin manifold with positive scalar $Scal.$ Suppose
$TM=T_1\oplus\cdots\oplus T_k, $ where $T_i$ are parallel distributions of dimension $n_i,i=1,\cdots,k,$ and $n_1> n_2\geq\cdots\geq n_k.$ Then one can define
a locally decomposable
Riemannian structure $J$ as follows
\beq
J\mid_{T_1}=Id,\ \  J\mid_{T_1^{\perp}}=-Id.
\eeq
Moreover,  we  define the following modified metric connection and corresponding
self-adjoint operator for $X\in \Gamma(TM^n)$ and $\phi\in\Gamma(\Sigma M^n)$
 \beq
\nabla^{\lambda_N}_X\phi
&:=&\nabla_X\phi+\frac{\lambda_N}{2n_1}X\cdot \phi+\frac{\lambda_N}{2n_1} J(X)\cdot \phi\\
\mathcal {T}^{\lambda_N}&:=& \nabla^{{\lambda_N}*}\nabla^{{\lambda_N}}+\frac{Scal}{4}Id.
\eeq

Compute for any eigenspinor $\psi_\alpha$ such that $D\psi_\alpha=\lambda_\alpha\psi_\alpha,\alpha=1,\cdots, N,$
\beQ
&{}&\int_{M^n}\langle\mathcal {T}^{\lambda_N}\psi_\alpha,\psi_\alpha\rangle\nonumber\\
&=&\left(\lambda_\alpha^2-\frac{2\lambda_\alpha\lambda_N}{n_1}+\frac{\lambda_N^2}{n_1}\right)
\int_{M^n}|\psi_\alpha|^2
+\frac{\lambda_N}{\lambda_\alpha}\cdot\frac{1}{2n_1}\int_{M^n}|D_J\psi_\alpha-D\psi_\alpha|^2.
\eeQ
But \cite{A1} gives
\beq
\|D\psi_\alpha\|^2
=\frac{n_1}{n_1-1}\|P\psi_\alpha\|^2+\frac{1}{n_1-1}\sum_{i=2}^k\varepsilon_i
\|D_i\psi_\alpha\|^2+\frac{n_1}{n_1-1}\int_{M^n}\frac{Scal}{4}|\psi_\alpha|^2\ \ \ \ \
\eeq
where $\varepsilon_i\triangleq\frac{n_1}{n_i}-1>0, \varepsilon_2\leq\cdots\leq \varepsilon_k,   i=2,\cdots,k$ and
$D_i$ is the ``partial" Dirac operator of subbundle $T_i.$
Hence
\beq\label{key inequality}
\|D\psi_\alpha\|^2
&\geq&\frac{1}{n_1-1}\sum_{i=2}^k\varepsilon_i
\|D_i\psi_\alpha\|^2+\frac{n_1}{4(n_1-1)}\min_M Scal\|\psi_\alpha\|^2.
\eeq
Suppose
$$\{e_1,\cdots,e_{n_1},e_{n_1+1}, \cdots,e_{n_1+n_2 },\cdots,e_{n_1+n_2+\cdots+n_{k-1}+1 },\cdots, e_n\}$$
is an adapted local orthnomal frame, i.e., such that
 $\{e_{n_1+n_2+\cdots+n_{i-1}+1},\cdots,e_{n_1+n_2+\cdots+n_i}\}$ spans $T_i.$
 Let $\alpha_i=n_1+n_2+\cdots+n_i,$
then
\beq
D_i&=&\sum_{\alpha={\alpha_{i-1}+1}}^{\alpha_{i}} e^\alpha\cdot\nabla_\alpha\\
D&=&\sum_{i=1}^kD_i\\
D_J&=&D_1-(D_2+\cdots+D_k)
\eeq
Note for $i\neq j,$ one has
\beq
D_iD_j+D_jD_i=0
\eeq
that is,
\beq
D^2=\sum_{i=1}^kD^2_i
\eeq
which, in turn, implies that
\beQ
\int_M|D_J\psi_\alpha-D\psi_\alpha|^2
&=&4\int_M\left|\sum_{i=2}^kD_i\psi_\alpha\right|^2\\
&=&4\sum_{i=2}^k\int_M\left|D_i\psi_\alpha\right|^2\\
&=&4\sum_{i=2}^k\varepsilon_i^{-1}\varepsilon_i\int_M\left|D_i\psi_\alpha\right|^2\\
&\leq&4\varepsilon_2^{-1}\sum_{i=2}^k\varepsilon_i\int_M\left|D_i\psi_\alpha\right|^2.
\eeQ
Hence this, together with the inequality (\ref{key inequality}) implies
\beQ
\|D_J\psi_\alpha-D\psi_\alpha\|^2\leq4\varepsilon_2^{-1}(n_1-1)
\left(\lambda_\alpha^2-\frac{n_1}{4(n_1-1)}{Scal_{min}}\right)\|\psi_\alpha\|^2.
\eeQ
This yields that
\beQ
&{}&\frac{\int_{M^n}\langle \mathcal {T}^{\lambda_N} \psi_\alpha,\psi_\alpha\rangle}{\int_{M^n}|\psi_\alpha|^2}\nonumber\\
&\leq&\left(\lambda_\alpha^2-\frac{2\lambda_\alpha\lambda_N}{n_1}+\frac{\lambda_N^2}{n_1}\right)
  +\frac{\lambda_N}{\lambda_\alpha}\cdot\frac{1}{2n_1}\cdot
  4\varepsilon_2^{-1}(n_1-1)
  \left(\lambda_\alpha^2-\frac{n_1}{4(n_1-1)}{Scal_{min}}\right)\\
&\leq&\left(1+2\varepsilon_2^{-1}\right)\frac{n_1-1}{n_1}\lambda_N^2-\frac{1}{2}
\varepsilon_2^{-1}Scal_{min}.
\eeQ

\begin{theorem}
Let $M$ be a compact Riemannian spin manifold with positive scalar $Scal.$ Let
$TM=T_1\oplus\cdots\oplus T_k, $ where $T_i$ are parallel distributions of dimension $n_i,i=1,\cdots,k,$ and $n_1> n_2\geq\cdots\geq n_k.$ Then for any positive integer N,
$$\lambda_N(D)^2\geq \frac{n_1-n_2}{n}\frac{n_1}{n_1-1}c\lambda_{k+1}(\Delta)+\frac{n_1}{4(n_1-1)}Scal_{\min},$$
where $k=\left[\frac{N}{2^{[\frac{n}{2}]}+1}\right],   c=\frac{1}{8(2^{[\frac{n}{2}]}+1)^2}.$
\end{theorem}

\emph{{\textbf{Proof}}}. As before,
the min-max principle yields
\beQ
c\lambda_{k+1}(\Delta)+\frac{1}{4}Scal_{\min}
&\leq&\lambda_N(\mathcal {T}^{\lambda_N})\\
&\leq&\max_{\alpha\in \{1,\cdots,N \}}
\frac{\int_{M^n}\langle\mathcal{T}^{\lambda_N}\psi_\alpha,\psi_\alpha\rangle}{\int_{M^n}|\psi_\alpha|^2}\\
&\leq&\frac{n}{n_1-n_2}\frac{n_1-1}{n_1}\lambda_N^2-
\frac{n_2}{2(n_1-n_2)}Scal_{min},
\eeQ
so we are done.

\end{document}